\documentstyle[twoside,amstex]{article}
\input amssym.def
\input amssym.tex
\def\bc{\begin{center}}
\def\ec{\end{center}}
\def\no{\noindent}

\begin{document}
\thispagestyle{empty} \vspace*{3 true cm} \pagestyle{myheadings}
\markboth {\hfill {\sl Huanyin Chen}\hfill} {\hfill{\sl Internal
cancellation over SSP rings}\hfill} \vspace*{-1.5 true cm}
\bc{\large\bf Internal Cancellation over SSP Rings}\ec

\vskip6mm
\bc{{\bf Huanyin Chen}\\[1mm]
Department of Mathematics, Hangzhou Normal University\\
Hangzhou 310036, People's Republic of China}\ec

\vskip4mm \begin{abstract} A ring is SSP if the sum of two direct
summands is a direct summand. A ring has internal cancellation if
every its (von Neumann) regular elements are unit-regular. We show
that in an SSP ring having internal cancellation, any regular
element is special clean. Our main results also imply that for any
SSP ring internal cancellation and idempotent stable range $1$
coincide with each other. Internal cancellation over SSP was then
characterized by special clean elements.


\end{abstract}

\vskip10mm\section{Introduction} \vskip4mm Let $R$ be a ring with
an identity. An element $a\in R$ is (unit) regular if there exists
some (unit) $x\in R$ such that $a=axa$. A ring $R$ is (unit)
regular if and only if every element in $R$ is (unit) regular. As
is well known, a ring $R$ is unit-regular if and only if every
element in $R$ is the product of an idempotent and a unit. For
general theory of regular rings, we refer the reader to
~\cite{GO}. An element $a\in R$ is (special) clean if it is the
sum of an idempotent $e$ ($aR\bigcap eR=0$) and a unit $u$. In
~\cite[Theorem 5]{CY}, Camillo and Yu claimed that every element
in unit-regular rings is clean. Unfortunately, there was a gap in
their proof. In ~\cite[Theorem 1]{CK}, Camillo and Khurana
improved ~\cite[Theorem 5]{CY} and proved that a ring $R$ is
unit-regular if and only if every element in $R$ is special clean.

A ring $R$ has stable range 1 provided that $Ra+Rb=R$ with $a,b\in
R$ implies that $a+zb\in R$ is a unit. A regular ring $R$ is
unit-regular if and only if $R$ has stable range 1. An interesting
problem is to extend Camillo-Khurana's theorem to certain rings
having stable range 1. A ring $R$ is an exchange ring if for any
$a\in R$ here exists an idempotent $e\in aR$ such that $1-e\in
(1-a)R$. Clearly, every regular ring is an exchange ring. In
~\cite[Theorem 2.1]{CH}, the author proved that an exchange ring
$R$ has stable range 1 if and only if every regular element in $R$
is special clean.

In ~\cite{KL}, Khurana and Lam introduced internal cancellation
(IC, for short) over a ring. A ring $R$ is IC if every regular
element in $R$ is unit-regular. Obviously, a regular ring $R$ is
unit-regular if and only if $R$ is IC. As is well known, a ring
$R$ is IC if and only if for any regular $a,b\in R$, $Ra+Rb=R$
implies that there exists a $z\in R$ such that $a+zb\in R$ is a
unit. Clearly, stable range 1 implies internal cancellation, but
the converse is not true, e.g. ${\Bbb Z}$. For an exchange ring
$R$, stable range 1 and internal cancellation coincide with each
other~\cite[Lemma 1.3.1]{CH1}. Many types of IC rings are
investigated in ~\cite{KL}. A further attractive problem is to
extend Camillo-Khurana's theorem to certain rings having IC. In
~\cite{WCK}, Wang et al. consider a kind of idempotent stable
range one for rings having IC. In ring theoretic version, the main
result in~\cite{WCK} stated

\vskip4mm \hspace{-1.8em} {\bf Theorem 1.1.}\ \ Let $R$ be a ring,
let $a\in R$ is regular, and let $end_R(r(a))$ or $end_R(aR)$ has
stable range 1. If $R$ is IC, then $Ra+Rb=R$ implies that there
exists an idempotent $e\in R$ such that $a+eb\in U(R)$ and
$R=aR\oplus eR$.

\vskip2mm \no Then Camillo-Khurana's theorem will be a special
case of this one. But stable range 1 is necessary in the preceding
theorem, as the substitution of modules is essential in its proof.
The motivation of this short note is to extend Camillo-Khurana's
theorem in a broader context with no restriction of stable range
1.

Following Garcia~\cite{G}, a ring $R$ is said to being the summand
sum property (briefly SSP) if the sum of two direct summands of
$R_R$ is also a direct summand of $R$. For instance, every (von
Neumann)regular ring and every every abelian ring (idempotents are
central, e.g., commutative rings, right duo rings) are SSP. We
shall prove that for any SSP ring internal cancellation and
idempotent stable range $1$ coincide with each other. Then in an
SSP ring having internal cancellation, any regular element is
special clean. Internal cancellation over SSP was then
characterized by special clean elements.

Throughout, all rings are associative with an identity. Let $R$ be
a ring. $U(R)$ will denote the set of all units in $R$. The right
annihilator $r(a)=\{~r\in R~|~ar=0\}$.

\section{The main results}

We begin with several lemmas which will be needed in our proof of
the main results.

\vskip4mm \hspace{-1.8em} {\bf Lemma 2.1.}\ \ {\it Every special
clean element in a ring is unit-regular.} \vskip2mm\hspace{-1.8em}
{\it Proof.}\ \ Let $a\in R$ be special clean. Then there exists
an idempotent $e\in R$ and a unit $u\in R$ such that $a=e+u$ and
$aR\bigcap eR=0.$ Hence, $au^{-1}=eu^{-1}+1$. Thus,
$au^{-1}e=eu^{-1}e+e\in aR\bigcap eR=0$. This yields
$au^{-1}(a-u)=0$, and so $au^{-1}a=a$. Therefore $a\in R$ is
unit-regular.\hfill$\Box$

\vskip4mm \hspace{-1.8em} {\bf Lemma 2.2.}~\cite[Proposition
2.2]{WCK} \ \ Let $M$ a right $k$-module, and let $R=End_k(M)$.
Then two submodules $A$ and $B$ of $M$ has a common direct
complement if and only if there exists an idempotent $e\in R$ such
that $eM=A$ and $e~|_B: B\to A$ is an isomorphism.

\vskip4mm \hspace{-1.8em} {\bf Lemma 2.3.}~\cite[Ehrlich¡¯s
Theorem 1.1]{KL} \ \ A ring $R$ is IC if and only if $R=A_1\oplus
B=A_2\oplus C$ with $A_1\cong A_2$ implies $B\cong C$.

\vskip4mm \hspace{-1.8em}
 {\bf Theorem 2.4.}\ \ {\it Let $R$ be
SSP. Then the following are equivalent:}\vspace{-.5mm}
\begin{enumerate}
\item [(1)] {\it $R$ is IC;}
\vspace{-.5mm} \item [(2)] {\it For any regular $a,b\in R$,
$Ra+Rb=R$ implies that there exists an idempotent $e\in R$ such
that $a+eb\in U(R)$ and $aR\oplus eR=R$.} \vspace{-.5mm}
\item [(3)] {\it Every regular element in $R$ is special
clean.}\end{enumerate} \vspace{-.5mm} {\it Proof.}\ \
$(1)\Rightarrow (2)$ Suppose that $Ra+Rb=R$ and $a,b\in R$ are
regular. Write $a=axa$ and $x=xax$ for some $x\in R$. Let $K=r(a),
D=xaR, I=aR$ and $C=(1-ax)R$. Then $R=K\oplus D=C\oplus I$ with
$a|_D:D\cong I$. Since $R$ is IC, we see that $K\cong C$. Clearly,
$K\bigcap r(b)=0$. Hence, $b: K\to bK$ is an isomorphism, and then
$C\cong K\cong bK$.

Since $Ra+Rb=R$, we see that $Rb(1-xa)=R(1-xa)$. We infer that
$b(1-xa)\in R$ is regular, and then $bK=b(1-xa)R$ is a direct
summand of $R$. Thus, we have an idempotent $f\in R$ such that
$bK=fR$. Set $g=1-xa$. Then $C=gR$ and $fR\cong gR$. As $R$ is IC,
we get $(1-f)R\cong (1-g)R$. In view of \cite[Proposition 2.5]{G},
$R$ has the summand intersection property, i.e., the intersection
of two direct summands of $R$ is a direct summand. Thus,
$fR\bigcap gR$ is a direct summand of $R$. Write $R=fR\bigcap
gR\oplus L$. Then $fR=fR\bigcap gR\oplus fR\bigcap L$ and
$gR=fR\bigcap gR\oplus gR\bigcap L$. This shows that
$$R=fR\bigcap gR\oplus (1-f)R\oplus fR\bigcap L=fR\bigcap gR\oplus
(1-g)R\oplus gR\bigcap L.$$ As $fR\bigcap gR\oplus (1-f)R\cong
fR\bigcap gR\oplus (1-g)R$, we get $\varphi: fR\bigcap L\cong
gR\bigcap L$.

Let $E=\{ x+\varphi(x)\mid x\in fR\bigcap L\}$. Clearly,
$fR+gR=fR+E$. If $x+\varphi(x)\in fR\bigcap E$ for $x\in fR\bigcap
L$, then $\varphi(x)=\big(x+\varphi(x)\big)-x\in \big(gR\bigcap
L\big)\bigcap fR=\big(fR\bigcap gR\big)\bigcap L=0$. As $\varphi$
is an isomorphism, we get $x=0$, and so $x+\varphi(x)=0$. We infer
that $fR\bigcap E=\emptyset$. Thus, $fR+gR=fR\oplus E$. Likewise,
$fR+gR=gR\oplus E$. Since $R$ is SSP, there exists a right
$R$-module $F$ such that $$R=(fR+gR)\oplus F=bK\oplus (E\oplus
F)=C\oplus (E\oplus F).$$ Hence, $C$ and $bK$ has a common direct
complement. In view of Lemma 2.2, we have an idempotent $e\in R$
such that $eR=C$ and $e: bK\to C$ is an isomorphism. Clearly,
$R=aR\oplus C=aR\oplus eR$. Let $u:=a+eb: R=K\oplus D\to R,
x+y\mapsto a(y)+eb(x+y)$ for any $x\in K,y\in D$.

One easily checks that the diagram of exact sequences
$$\begin{array}{ccccccccc}
0&\to &K&\to &K\oplus D&\to &D&\to &0\\
&&eb \downarrow& &a+eb\downarrow& &a\downarrow& &\\
0&\to &C&\to &C\oplus I&\to &I&\to &0.
\end{array}$$ commutes. As $eb: K\to C$ and $a: D\to I$ are both
isomorphisms, by using the Five Lemma, we conclude that $a+eb\in
U(R)$, as asserted.

$(2)\Rightarrow (3)$ For any regular $a\in R$, we see that
$Ra+R(-1)=R$. Thus, we can find an idempotent $e\in R$ such tht
$a-e\in U(R)$ and $aR\bigcap eR=0$. That is, $a\in R$ is special
clean.

$(3)\Rightarrow (1)$ In view of Lemma 2.1, every regular element
in $R$ is unit-regular. Therefore $R$ is IC, as
asserted.\hfill$\Box$

\vskip4mm \hspace{-1.8em} {\bf Remark 2.5.}\ \ From a similar
proof of Theorem 2.4, we actually derive a stronger version: An
SSP ring $R$ is IC if and only if for any regular $a,b\in R$,
$r(a)\bigcap r(b)=0$ implies that there exists an idempotent $e\in
R$ such that $a+eb\in U(R)$ and $aR\oplus eR=R$. As SSP and IC are
both left-right symmetric properties, we see that an SSP ring $R$
is IC if and only if for any regular $a,b\in R$, $aR+bR=R$ implies
that there exists an idempotent $e\in R$ such that $a+be\in U(R)$
and $Ra\oplus Re=R$, by applying Theorem 2.4 to the opposite ring
$R^{op}$ of $R$.

\vskip4mm \hspace{-1.8em} {\bf Corollary 2.6} [Camillo-Khurana's
theorem]\ \ {\it Let $R$ be a ring. Then the following are
equivalent:}\vspace{-.5mm}
\begin{enumerate}
\item [(1)] {\it $R$ is unit-regular.} \vspace{-.5mm} \item [(2)] {\it
Every element in $R$ is special clean.}\end{enumerate}
\vspace{-.5mm} {\it Proof.}\ \ This is obvious by Theorem 2.4 and
Lemma 2.1.\hfill$\Box$

\vskip4mm \hspace{-1.8em} {\bf Corollary 2.7.}\ \ {\it Let $R$ be
abelian. Then for any regular $a\in R$, there exists a unique
idempotent $e\in R$ and a unit $u\in R$ such that $a=e+u$ and
$aR\bigcap eR=0.$} \vskip2mm\hspace{-1.8em} {\it Proof.}\ \
Clearly, $R$ is SSP and IC. Let $a\in R$ be regular. Then $a\in R$
is special clean, by Theorem 2.4. Assume that there exist
idempotents $e,e'\in R$ and units $u,u'\in R$ such that
$a=e+u=e'+u'$ and $aR\bigcap eR=0=aR\bigcap e'R.$ Then
$au^{-1}=eu^{-1}+1$. As $R$ is abelian, we see that
$au^{-1}(1-e)=1-e$, and so $au^{-1}(1-e)e'=e'(1-e)\in aR\bigcap
e'R=0$. This forces that $e'=e'e$. Likewise, $e=ee'$. Therefore
$e=e'$, as needed.

\vskip4mm \hspace{-1.8em} {\bf Example 2.8.}\ \ The ring ${\Bbb
Z}$ of all integers is an SSP ring having IC, but it has not
stable range 1.

\vskip4mm \hspace{-1.8em} {\bf Theorem 2.9} \ \ {\it Let $R$ be a
ring. Then the following are equivalent:}\vspace{-.5mm}
\begin{enumerate}
\item [(1)] {\it $R$ is SSP and IC.}
\vspace{-.5mm}
\item [(2)] {\it
The product of two regular elements is unit-regular.}
\vspace{-.5mm}
\item [(3)] {\it
The product of two regular elements is special
clean.}\end{enumerate} \vspace{-.5mm} {\it Proof.}\ \
$(1)\Rightarrow (3)$ Let $a,b\in R$ be regular. Since $R$ is IC,
we can find idempotents $e,f\in R$ and units $u,v\in R$ such that
$a=ue$ and $b=fv$. Then $ab=u(ef)v$. Observing that
$$(1-e)R\oplus efR=(1-e)R+fR=gR~\mbox{for some idempotent}~g\in
R,$$ we get $R=efR\oplus (1-e)R\oplus (1-g)R$. This implies that
$ef\in R$ is regular. Hence, $u(ef)v$ is regular, and then so is
$ab$ in $R$. By virtue of Theorem 2.4, $ab\in R$ is special clean,
as required.

$(3)\Rightarrow (2)$ This is clear in terms if Lemma 2.1.

$(2)\Rightarrow (1)$ Clearly, every regular element in $R$ is
unit-regular. That is, $R$ is IC. Let $e,f\in R$ be idempotents.
Then $eR+fR=eR+(1-e)fR$. By hypothesis, there exists an idempotent
$g\in R$ such that $(1-e)fR=gR$. Then we check that $eR+fR=(e+g)R$
where $e+g\in R$ is an idempotent. Therefore $eR+fR$ is a direct
summand of $R$, hence that $R$ is SSP.\hfill$\Box$

\vskip4mm \hspace{-1.8em} {\bf Corollary 2.10} \ \ {\it Let $R$ be
a ring. Then the following are equivalent:}\vspace{-.5mm}
\begin{enumerate}
\item [(1)] {\it $R$ is SSP and IC.}
\vspace{-.5mm}
\item [(2)] {\it
The product of finitely many of regular elements is unit-regular.}
\vspace{-.5mm}
\item [(3)] {\it
The product of finitely many elements is special
clean.}\end{enumerate} \vspace{-.5mm} {\it Proof.}\ \
$(1)\Rightarrow (3)$ Let $a=a_1\cdots a_{n-1}a_n (n\geq 2)$ where
each $a_i\in R$ is regular. Then $a_1a_2\in R$ is regular, by
Theorem 2.9. By iteration of this process, we see that
$a_1a_2\cdots a_{n-1}\in R$ is regular. In light of Theorem 2.9,
we conclude that $(a_1\cdots a_{n-1})a_n\in R$ is special clean.
Therefore $a\in R$ is special clean, as desired.

$(3)\Rightarrow (2)$ This is obvious by Lemma 2.1.

$(2)\Rightarrow (1)$ is trivial from Theorem 2.9.\hfill$\Box$

\vskip4mm We note that the number of regular elements in Theorem
2.9 can not be reduced from two to one as the following shows.

\vskip4mm \hspace{-1.8em} {\bf Example 2.11.}\ \ Let $R=\left(
\begin{array}{cc}
{\Bbb Z}_3&{\Bbb Z}_3\\
0&{\Bbb Z_3}
\end{array}
\right)$. One easily checks that $R$ has stable range one. In view
of Theorem 1.1, every regular element in $R$ is special clean.
Choose $e=\left(
\begin{array}{cc}
1&1\\
0&0
\end{array}
\right)$ and $f=\left(
\begin{array}{cc}
0&1\\
0&1
\end{array}
\right)$. Then $e,f\in R$ are idempotents, and then regular. But
$ef=\left(
\begin{array}{cc}
0&-1\\
0&0
\end{array}
\right)\in R$ is not regular, as $ef\not \in (ef)R(ef)=0$. Thus,
$ef\in R$ is not special clean.

\end{document}